\input amstex
\input epsf
\documentstyle{amsppt}
\NoBlackBoxes
\nologo
\tolerance 1000

\define\CP{\Bbb C\bold P}
\define\Aut{\operatorname{Aut}}
\define\st{\operatorname{st}}
\define\MO{\overline{\Cal M}}
 
\define \R{\Bbb R} 

\define\gs{\geqslant}

\define\T{\Cal T}
\define\C{\Bbb C}
\define\I{\Cal I}
\define\E{\Cal E}
\define\Hu{\Cal H}
\define\bHu{\overline{\Cal H}}

\define\exsum{\mathop{{\sum}^\varnothing}}
\def\ll{{\ell\negthickspace\ell}}
\def\bll{{\ll}}

\topmatter

\title
On double Hurwitz numbers in genus $0$
\endtitle

\author
S.~Shadrin$^*$, M.~Shapiro$^{**}$ and A.~Vainshtein$^{\dag}$
\endauthor

\affil
$^*$ Department of Mathematics, University of Zurich\\
Winterthurerstrasse 190, CH-8057 Zurich, Switzerland\\
{\tt sergey.shadrin\@math.unizh.ch}\\
and Department of Mathematics, Institute of System Research\\
Nakhimovskii prosp.~36-1, 117218 Moscow, Russia\\
$^{**}$ Department of Mathematics, Michigan State University\\
East Lansing, MI 48824, USA, {\tt mshapiro\@math.msu.edu}\\
$^\dag$ Depts. of Mathematics and Computer Science, University of Haifa
\\ Mount Carmel, 31905 Haifa, Israel, {\tt alek\@cs.haifa.ac.il}
\endaffil

\abstract
We study double Hurwitz numbers in genus zero 
counting the number of covers $\CP^1\to\CP^1$ with two 
branching points with a given branching behavior. By the recent result
due to Goulden, Jackson and Vakil, these numbers are 
piecewise polynomials in the multiplicities of the preimages of the 
branching points. We describe the partition of the parameter space into
polynomiality domains, called chambers, and provide an expression for the 
difference of two such polynomials for two neighboring chambers. Besides,
 we provide an explicit formula for the polynomial in a certain
chamber called totally negative, which enables us to calculate double 
Hurwitz  numbers 
in any given chamber as the polynomial for the totally negative
chamber plus the
sum of the differences between the neighboring polynomials along a path 
connecting the totally negative chamber with the given one.  
\endabstract
\rightheadtext{Double Hurwitz numbers}
\leftheadtext{S.~Shadrin, M.~Shapiro and A.~Vainshtein}
\endtopmatter

\document
\heading {\S 1. Introduction and results} \endheading

Let $\mu=(\mu_1,\dots,\mu_m)$ be a partition of an integer $d$.
The {\it Hurwitz number\/}  $h^g_{\mu}$
is the number of genus $g$ branched covers of $\CP^1$ with branching
corresponding to $\mu$  over a fixed point (usually 
identified with $\infty$) and an appropriate number
of fixed simple branching points. 
Recall that each cover is counted with the weight $c^{-1}$, 
where $c$ is the number of automorphisms of the cover. 

Hurwitz numbers possess a rich structure explored by many authors in different
fields, such as algebraic geometry, representation theory, integrable 
systems, combinatorics,
and mathematical physics. Being very far from trying to describe all these
achievements, we mention only the so-called ELSV-formula (see \cite{ELSV1,
ELSV2, GrVa})
$$
h^g_{\mu}=\frac{(m+d+2g-2)!}{|\Aut \mu|}\prod_{i=i}^m\frac{\mu_i^{\mu_i}}
{\mu_i!}
\int_{\MO_{g,m}}\frac{c(\Lambda^\vee)}
{(1-\mu_1\psi_1)\cdots(1-\mu_m\psi_m)},
$$
where  $\Aut\mu$ is the group of symmetries of the set $\mu$, $\MO_{g,m}$ is 
the Deligne--Mumford compactification of the moduli space of genus $g$
curves with $m$ marked points, $\psi_i$ is the first Chern class of
the cotangent bundle at the $i$th marked point, and
$c(\Lambda^\vee)$ is the total Chern class of the dual to the Hodge bundle
(see the above mentioned papers for exact definitions).

Let $\nu=(\nu_1,\dots,\nu_n)$ be another partition
of the same integer $d$. The {\it double Hurwitz number\/} $h^g_{\mu;\nu}$
is the number of genus $g$ branched covers of $\CP^1$ with branchings
corresponding to $\mu$ and $\nu$ over two fixed points (in what follows we
identify them with $\infty$ and $0$, respectively) and an appropriate number
of fixed simple branching points. We denote the latter number $r_{\mu;\nu}^g$; 
by the Riemann--Hurwitz theorem, $r_{\mu;\nu}^g=m+n+2g-2$. To simplify the 
exposition, we assume that the points mapped to $\infty$ and $0$ are
labelled, so the double Hurwitz numbers under this convention are 
$|\Aut \mu||\Aut \nu|$ larger than they would be under the usual convention. 
 

Most of the known results concerning double Hurwitz numbers treat only the 
so-called one-part (or polynomial) case, when $m=1$ and $n$ is arbitrary.  
One-part double Hurwitz numbers in genus zero where studied in \cite{ZL}
(see also \cite{Z1} for an earlier version of the same result). It is proved 
there that 
$$
h^0_{(d);\nu}=(n-1)!d^{n-2}; \tag 1.1
$$
recall that $n-1=r^0_{(d);\nu}$. In fact, a combinatorial result much more
general than (1.1) was obtained already in \cite{GJ} using Lagrange inversion; 
in \cite{ZL} the same formula was reproved
by methods of algebraic geometry, and in this way its algebro-geometric
meaning was clarified.

Another way to get (1.1) is suggested in
\cite{Sh1}. The aim of the paper is to derive an ELSV-type formula for 
$h^g_{(d);\nu}$. Two such formulas are derived; the first of them
(see \cite{Sh1, Sh3}) represents $h^g_{(d);\nu}$ as 
$$
h^g_{(d);\nu}=r!d^{r-1}\int_{\Delta}\psi_1^{r-1},
$$
where $r=r^g_{(d);\nu}=n+2g-1$, 
$\Delta$ is  a certain cycle in the moduli space $\MO_{g,n+1}$ and
$\psi_1$ has the same meaning as in the ELSV-formula.
In the case $g=0$ the cycle $\Delta$
coincides with the whole space $\MO_{0,n+1}$, and since the integral is known
to be equal to~$1$, we immediately get (1.1). The second formula 
(see \cite{Sh1, Sh2}) represents $h^0_{(d);\nu}$ as 
$$
h^0_{(d);\nu}=d^{n-2}\int_{\MO_{0,d(n-1)+2}}\psi_1^{n-2}\Psi_2\cdots\Psi_{n+1}.
$$
The definition of the classes $\Psi_i$ is recursive and rather involved. 
The only example considered in \cite{Sh1} is $n=2$, in which case, after
some effort, one gets the correct answer. 

A different ELSV-type formula for the one-part double Hurwitz numbers
is conjectured in \cite{GJVa}. It is proved there for the case $g=0$. 
The answer look as follows:
$$
h^0_{(d);\nu}={(n-1)!}d\int_{\MO_{0,n}}
\frac1{(1-\nu_1\psi_1)\cdots(1-\nu_n\psi_n)},
$$
where $\psi_i$ are the same as in the ELSV formula. For higher
genera, the conjectured formula looks very similar, however, it is yet not
clear over which moduli space should the integral be taken. 
On the computational side, \cite{GJVa} suggests the following result,
which is a common generalization of (1.1) and \cite{SSV, Th.~6}:
$$
h^g_{(d);\nu}={r!}d^{r-1}[t^{2g}]\prod_{k\gs1}
\left(\frac{\sinh{kt/2}}{kt/2}\right)^{c_k},
$$
where $r=r^g_{(d);\nu}$, $\{c_k\}$ is a finite sequence defined solely
by $\nu$, and $[A]B$ is the coefficient of $A$ in $B$. This result extends 
previous computations carried out in \cite{Ku} for $n=1,2$.

Much less is known for double Hurwitz numbers with arbitrary $\mu$ and $\nu$. 
In fact, there are only two
general results. First, it is proved in \cite{Ok} that the exponent of the
generating function for the double Hurwitz numbers is a $\tau$-function
for the Toda hierarchy of Ueno and Takasaki. Second,  Theorem~2.1 in
\cite{GJVa} states that for fixed $g$, $m$, and $n$, double 
Hurwitz numbers
are piecewise polynomial in variables $\mu_1,\dots,\mu_m,\nu_1,\dots,
\nu_n$, and that the highest degree of this piecewise polynomial is 
constant and equal to $m+n+4g-3$. Moreover, it is proved in \cite{GJVa} that 
for genus zero case this piecewise polynomial is homogeneous.

In a different direction, Theorem~4.1 of the same paper treats
$h^g_{\mu;\nu}$ as a function of $g$ for fixed partitions $\mu$ and $\nu$.
In principle, any such function can be obtained by recursive computation;
the only case done explicitly is $m=n=2$ and $\mu_1>\nu_1>\nu_2>\mu_2$.
Similar results are also obtained in \cite{Ku} for all double Hurwitz numbers
of degree at most $5$.

Still another approach suggested in \cite{GJVa} provides exact formulas
for double Hurwitz numbers in genus zero when $m=2$ or $m=3$. In particular,
$$
h^0_{\mu_1,\mu_2;\nu}=\frac{n!|\Aut\nu|}{d}\sum
\frac{l(\rho)!\prod_{j\gs1}\rho_j}{|\Aut\rho||\Aut\sigma||\Aut\tau|}
(\mu_1-|\sigma|)(\mu_2-|\tau|)\mu_1^{l(\sigma)-1}\mu_2^{l(\tau)-1},
\tag 1.2
$$
where
the summation is over partitions $\rho=(\rho_1,\dots,\rho_{l(\rho)})$,
$\sigma=(\sigma_1,\dots,\sigma_{l(\sigma)})$, 
$\tau=(\tau_1,\dots,\tau_{l(\tau)})$ with $\rho\cup\sigma\cup\tau=\nu$
and $|\sigma|=\sum_{i=1}^{l(\sigma)}\sigma_i<\mu_1$,
$|\tau|=\sum_{i=1}^{l(\tau)}\tau_i<\mu_2$, see \cite{GJVa, Corollary 5.11}.

Finally, the approach developed in \cite{Z2, KL} allows in principle to obtain 
formulas for $h^0_{\mu;\nu}$ as functions of the degree $d$ for a very 
special choice of the partitions $\mu$ and $\nu$: $\mu=(\alpha, 1,\dots,1)$,
$\nu=(\beta, 1,\dots,1)$, $\alpha$ and $\beta$ being fixed
partitions. For example,
for $\alpha=\beta=(3)$ one gets
$$
h^0_{3^11^{d-3},3^11^{d-3}}=\frac34(27d^2-137d+180)
\frac{d^{d-6}(2d-6)!}{(d-3)!}.
$$


In this note we study the homogeneous piecewise polynomial of degree $m+n-3$
defining double Hurwitz numbers $h^0_{\mu;\nu}$.
To formulate the results we need to introduce some 
notation. We consider $x_1=\mu_1,\dots,x_m=\mu_m,y_1=\nu_1,\dots,y_n=\nu_n$ 
as coordinates of a point in $\R^{m+n}$. 
The {\it parameter space\/} is the cone in  $\R^{m+n}$ given by the 
inequalities $x_1\ge
\dots\ge x_m\ge 0$,  $y_1\ge\dots\ge y_n\ge 0$ and the equality 
$\sum_{i=1}^mx_i=\sum_{j=1}^ny_j$. A {\it resonance\/} is a hyperplane
$x_I=y_J$, where $I\subset[1,m]$ and $J\subset[1,n]$ are proper subsets
and $a_K$ stands for $\sum_{i\in K}a_i$ for any sequence $a_1,\dots,a_k$ 
and any subset $K\subseteq[1,k]$. The connected components of the complement
to the union of all resonances are called {\it chambers}. 

\proclaim{Theorem 1.1} Let $(\mu,\nu)$ vary within a closure of a chamber, then
double Hurwitz numbers $h^0_{\mu;\nu}$ are given by a homogeneous 
polynomial of degree $m+n-3$.
\endproclaim

Consider a chamber $C$, and let $P_C$ be the corresponding homogeneous 
polynomial. The most convenient way to identify $C$ is to pick up a
reference point
$(\alpha,\beta)\in C$; in this case we write $P_{\alpha;\beta}$ instead of
$P_C$. Observe that the only role of the reference point $(\alpha,\beta)$ is to
indicate the choice of the chamber, and that its coordinates are not
necessary integers. It is clear that $P_{\mu;\nu}(\mu,\nu)=h^0_{\mu;\nu}$;
however, $P_{C}(\mu',\nu')$ may differ from $h^0_{\mu';\nu'}$
if $(\mu',\nu')\notin C$. 

The total number of resonances is
equal to $2(2^{m-1}-1)(2^{n-1}-1)$, since each resonance $x_I=y_J$ can be
also written as $x_{\bar I}=y_{\bar J}$, where  bar stands for the complement. 
In what follows we always assume that $1\notin I$. 
Therefore, each chamber is defined by a sequence
of $w=2(2^{m-1}-1)(2^{n-1}-1)$ signs of the expressions $x_I-y_J$; observe that
the total number of chambers is less than $2^w$, since certain combinations
of signs are impossible. We say that two chambers 
are {\it neighboring\/} along the resonance  $x_I=y_J$ if the corresponding 
sign sequences differ only in the position corresponding to this resonance.

Let $C$ and $C'$ be two chambers neighboring along the resonance  $x_I=y_J$; 
without loss of generality we assume that $x_I-y_J>0$ in $C$.

\proclaim{Theorem 1.2} Let $(\mu,\nu)$ be an arbitrary point
in $C$. Then
$$
P_C-P_{C'}=\binom{m+n-2}{|I|+|J|-1}(x_I-y_J)
P_{\mu(\bar I),\mu_I-\nu_J;\nu(\bar J)}
P_{\mu(I);\nu(J),\mu_I-\nu_J},
$$
where $a(K)$ is the subsequence of $a_1,\dots,a_k$ consisting of terms
$a_i$, $i\in K$.
\endproclaim

\remark{Remark} Here and in what follows we omit the arguments of 
polynomials whenever this does not lead to a confusion. Usually, the
arguments are formed from the components of $x$ and $y$ according to
the same rules as the coordinates of the reference point are formed
from the parts of $\mu$ and $\nu$. For example,
$P_{\mu(I);\nu(J),\mu_I-\nu_J}$ has $|I|+|J|+1$ arguments. The first
$|I|$ of them are $x_i$, $i\in I$, then follow $y_j$, $j\in J$, and
the last argument is $x_I-y_J$.
\endremark

Consider the {\it totally negative\/} chamber, that is, the one for which
all the signs in the corresponding sequence are negative. The following result
was conjectured in \cite{GJVa, Conjecture 5.10}.

\proclaim{Theorem 1.3} The polynomial corresponding to the totally negative
cell is given by $(m+n-2)!x_1^{n-1}(x_1+\dots+x_m)^{m-2}$.
\endproclaim

Theorems 1.1--1.3 give rise to recurrence relations expressing double
Hurwitz numbers of degree $d$ via double Hurwitz numbers of lesser
degrees. In general, these relations are rather cumbersome, however,
for $m=2$ one gets a very simple explicit formula, which is 
easier than~(1.2). The details of the corresponding computations are
presented in Section~2.
 
The authors are grateful to Max-Planck-Institut f\"ur Matematik, Bonn,
and to Institut Mittag-Leffler, Djursholm, Sweden, 
for hospitality in Summer 2004, and
Fall 2006. Our sincere thanks go to M.~Kazaryan, R.~Kulkarni and D.~Zvonkine
for useful discussions. S.S. was supported by the grants RFBR-05-01-01012-a,
RFBR-05-01-02806-CNRS-a, NSh-1972.2003.1, MK-5396.2006.1,
NWO-RFBR-047.011.2004.026 (RFBR-05-02-89000-NWO-a),
by the G{\"o}ran Gustafsson foundation, and by
Pierre Deligne's fund based on his 2004 Balzan prize in mathematics.
M.S. was supported by the grants DMS-0401178 and PHY-0555346. 
M.S and A.V. were supported by the grant BSF-2002375.

\heading {\S 2. Computations} \endheading

\subheading{2.1. General recurrence}
To find a recurrence relation for the double Hurwitz number $h^0_{\mu;\nu}$ via
Theorems~1.1--1.3, one has to pick a path connecting the totally
negative chamber with the chamber containing the point $(\mu,\nu)$ in
the parameter space. If $(\mu,\nu)$ lies on one or more resonances,
one can choose any of the adjacent chambers in an arbitrary
way. By multiplying all coordinates by a sufficiently big integer
$t>d^{m+n}$ and
perturbing slightly the resulting point within the same chamber, one
can ensure that the obtained point $(t\mu,\nu')$ is in general position, that
is, $\nu'_I=\nu'_J$ if and only if $I=J$. Pick the point
$(td-m(m-1)/2,m-1,m-2,\dots,1,\nu')$ in the totally negative chamber
and connect
it with $(t\mu,\nu')$ by the following path consisting of
$m-1$ segments. The first segment is of the form
$(td-m(m-1)/2-s,m-1+s,m-2,\dots,1,\nu')$, $s=0,1,\dots, t\mu_2-m+1$, the second
segment is of the form
$(td-t\mu_2-(m-1)(m-2)/2-s,t\mu_2,m-2+s,m-3,\dots,1,\nu')$,
$s=0,1,\dots, t\mu_3-m+2$, and so on. It is easy to see that 
each point on this path belongs to at most one resonance.
To formulate the recurrence relation, pick arbitrary numbers
$\varepsilon_3,\varepsilon_4,\dots,\varepsilon_m$ and $\delta$
satisfying inequalities
$0<\delta<\varepsilon_m<\dots<\varepsilon_3<1/m$ and denote
$\varepsilon=(\varepsilon_3,\dots,\varepsilon_m)$; 
the exact values
of  $\varepsilon_i$ and $\delta$ do not have any meaning, since these
numbers will be only used to indicate the corresponding chamber. Clearly, a
resonance $x_I=y_J$ is intersected by the above path, and hence
contributes to  $h^0_{\mu;\nu}$, if and only if $\mu_I>\nu_J$. If this
is the case, we consider $I=\{i_1,\dots,i_{|I|}\}$ and define
$k=\min\{j\in [1, |I|]: \mu_{i_1}+\dots+\mu_{i_j}>\nu_J\}$. Any $i\in
[1,m]$, $i\ne 1, i_k$, can be related to one of the four subsets:
$I_1=\{i\in I: i<i_k\}$, $\bar I_1=\{i\notin I: 1<i<i_k\}$,
$I_2=\{i\in I: i>i_k\}$,  $\bar I_2=\{i\notin I: i>i_k\}$. We thus get
the following result.

\proclaim{Theorem 2.1} Double Hurwitz numbers are given by
$$
\multline
h^0_{\mu;\nu}=(m+n-2)!d^{m-2}\mu_1^{n-1}+\sum_{\mu_I>\nu_J}
\binom{m+n-2}{|I|+|J|-1}(\mu_I-\nu_J)\times\\
P_{\xi_1,\mu(\bar I_1),\varepsilon(\bar I_2),\delta;\nu(\bar J)}
(\mu_1,\mu(\bar I_1\cup \bar I_2),\mu_I-\nu_J,\nu(\bar J))\times\\
P_{\mu(I_1),\xi_2,\varepsilon(I_2);\nu(J),\delta}
(\mu(I_1),\mu_{i_k}, \mu(I_2),\nu(J),\mu_I-\nu_J),
\endmultline
$$
where $\xi_1=\nu_{\bar J}-\mu_{\bar I_1}-\varepsilon_{\bar I_2}-\delta$,
$\xi_2=\nu_J-\mu_{I_1}-\varepsilon_{I_2}+\delta$.
\endproclaim

\subheading{2.2. Two-part double Hurwitz numbers} The general
expression in Theorem~2.1 looks very cumbersome. However, in the case of
two-part double Hurwitz numbers, when $m=2$, it can be written in a
very simple way. Indeed, in this case all the resonances are of the
form $x_2=y_J$, therefore $I=\{2\}$, $i_k=2$, $I_1=\bar I_1=I_2=\bar
I_2=\varnothing$. Therefore Theorem~2.1 yields
$$
\multline
h^0_{\mu_1,\mu_2;\nu}=n!\mu_1^{n-1}+\sum_{\mu_2>\nu_J}
\binom{n}{|J|}(\mu_2-\nu_J)\times\\
P_{\nu_{\bar J}-\delta,\delta;\nu(\bar J)}(\mu_1,\mu_2-\nu_J,\nu(\bar J))
P_{\nu_J+\delta;\nu(J),\delta}(\mu_{2},\nu(J),\mu_2-\nu_J).
\endmultline
$$
The second
polynomial in the right hand side corresponds to one-part double
Hurwitz numbers of total degree $\mu_2$ with $|J|+1$ zeros, hence, by~(1.1),
$P_{\nu_J+\delta;\nu(J),\delta}(\mu_{2},\nu(J),\mu_2-\nu_J)=|J|!\mu_2^{|J|-1}$.
The first polynomial in the right hand side corresponds to the totally
negative chamber for two-part double
Hurwitz numbers of total degree $d-\nu_J$ with $n-|J|$ zeros, hence, by
Theorem~1.3, $P_{\nu_{\bar J}-\delta,\delta;\nu(\bar J)}
(\mu_1,\mu_2-\nu_J,\nu(\bar J))=(n-|J|)!\mu_1^{n-|J|-1}$. 
Observe that the first summand in the
right hand side of the above formula can be also included in the
regular part of the sum for $J=\varnothing$. In what follows we
indicate this by writing $\exsum$ instead of $\sum$.

We thus obtain the following explicit formula for the two-part
double Hurwitz numbers in genus zero, which is simpler than~(1.2).
  
\proclaim{Corollary 2.2} The two-part double Hurwitz numbers are given by
$$
h^0_{\mu_1,\mu_2;\nu}=n!{\exsum_{\mu_2>\nu_J}}(\mu_2-\nu_J)
\mu_1^{n-|J|-1}\mu_2^{|J|-1}.
$$
\endproclaim

\subheading{2.3. Three-part double Hurwitz numbers} Consider now the case of
three-part double Hurwitz numbers, when $m=3$. Then the first segment of
the path intersects resonances of the form  $x_2=y_J$ and
$x_2+x_3=y_J$, while the second segment of the path intersects
resonances of the form  $x_3=y_J$ and $x_2+x_3=y_J$. Therefore, 
we have the following four types of intersection.

Type 1. $I=\{2\}$, $i_k=2$, $I_1=\bar I_1=I_2=\varnothing$, $\bar
   I_2=\{3\}$.

By Theorem 2.1, the contribution of such an intersection equals
$$
\multline
\binom{n+1}{|J|}(\mu_2-\nu_J)
P_{\nu_{\bar J}-\varepsilon_3-\delta,\varepsilon_3,\delta;\nu(\bar J)}
(\mu_1,\mu_3,\mu_2-\nu_J,\nu(\bar J))\times\\
P_{\nu_J+\delta;\nu(J),\delta}(\mu_{2},\nu(J),\mu_2-\nu_J).
\endmultline
$$
The second polynomial in the above expression is the same as in the
case of two-part double Hurwitz numbers, and its value is equal to 
$|J|!\mu_2^{|J|-1}$.
The first polynomial corresponds to the totally
negative chamber for three-part double
Hurwitz numbers of total degree $d-\nu_J$ with $n-|J|$ zeros, hence, by
Theorem~1.3, $P_{\nu_{\bar
J}-\varepsilon_3-\delta,\varepsilon_3,\delta;\nu(\bar J)}
(\mu_1,\mu_3,\mu_2-\nu_J,\nu(\bar J))=(n-|J|+1)!\mu_1^{n-|J|-1}(d-\nu_J)$.
Therefore, the total contribution of all intersections of type 1
equals
$(n+1)!\sum_{\mu_2>\nu_J}(\mu_2-\nu_J)\mu_1^{n-|J|-1}\mu_2^{|J|-1}(d-\nu_J)$. 

Type 2. $I=\{2,3\}$, $i_k=2$, $I_1=\bar I_1=\bar I_2=\varnothing$, $I_2=\{3\}$.

By Theorem 2.1, the contribution of such an intersection equals
$$
\multline
\binom{n+1}{|J|+1}(\mu_2+\mu_3-\nu_J)
P_{\nu_{\bar J}-\delta,\delta;\nu(\bar J)}
(\mu_1,\mu_2+\mu_3-\nu_J,\nu(\bar J))\times\\
P_{\nu_J-\varepsilon_3+\delta,\varepsilon_3;\nu(J),\delta}
(\mu_{2},\mu_3,\nu(J),\mu_2+\mu_3-\nu_J).
\endmultline
$$
The first polynomial in the above expression we have already
encountered in the case of two-part double Hurwitz numbers, and its
value is equal to  $(n-|J|)!\mu_1^{n-|J|-1}$.
The second polynomial corresponds to the chamber neighboring with the totally
negative chamber for two-part double
Hurwitz numbers of total degree $\mu_2+\mu_3$ with $|J|+1$ zeros,
hence, by Corollary~2.2,
$$
\multline
P_{\nu_J-\varepsilon_3+\delta,\varepsilon_3;\nu(J),\delta}
(\mu_{2},\mu_3,\nu(J),\mu_2+\mu_3-\nu_J)=\\
(|J|+1)!\left(\mu_2^{|J|}+ 
(\mu_3-(\mu_2+\mu_3-\nu_J))\mu_2^{|J|-1}\right)=(|J|+1)!\mu_2^{|J|-1}\nu_J.
\endmultline
$$
Therefore, the total contribution of all intersections of type 2
equals
$(n+1)!\sum_{\mu_2>\nu_J}(\mu_2+\mu_3-\nu_J)\mu_1^{n-|J|-1}\mu_2^{|J|-1}\nu_J$.
 
Type 3. $I=\{3\}$, $i_k=3$, $I_1=I_2=\bar I_2=\varnothing$, $\bar I_1=\{2\}$.

By Theorem 2.1, the contribution of such an intersection equals
$$
\multline
\binom{n+1}{|J|}(\mu_3-\nu_J)
P_{\nu_{\bar J}-\mu_2-\delta,\mu_2,\delta;\nu(\bar J)}
(\mu_1,\mu_2,\mu_3-\nu_J,\nu(\bar J))\times\\
P_{\nu_J+\delta;\nu(J),\delta}(\mu_3,\nu(J),\mu_3-\nu_J).
\endmultline
$$
The second polynomial in the above expression is again the same as in
the case of two-part double Hurwitz numbers, and its value is equal to
$|J|!\mu_3^{|J|-1}$. The first polynomial corresponds to three-part
double Hurwitz numbers of total degree $d-\nu_J$ with $n-|J|$ zeros, 
computed in some
chamber intersected by the first segment of our path. To compute these
numbers we have to take into account only intersections of types 1 and 2.
By the above reasoning, we get
$$
\multline
P_{\nu_{\bar J}-\mu_2-\delta,\mu_2,\delta;\nu(\bar J)}
(\mu_1,\mu_2,\mu_3-\nu_J,\nu(\bar J))=\\(n-|J|+1)!
\mathop{\exsum_{K\subseteq \bar J}}_{\mu_2>\nu_K}(\mu_2-\nu_K)(d-\nu_J-\nu_K)
\mu_1^{n-|J|-|K|-1}\mu_2^{|K|-1}+\\
\mathop{\exsum_{K\subseteq \bar J}}_{\mu_2>\nu_K}(\mu_2+(\mu_3-\nu_J)-\nu_K)
\mu_1^{n-|J|-|K|-1}\mu_2^{|K|-1}\nu_K.
\endmultline
$$
Therefore, the total contribution of all intersections of type 3 equals
$$
(n+1)!\sum_{\mu_3>\nu_J}(\mu_3-\nu_J)\mu_3^{|J|-1}
\mathop{\exsum_{K\subseteq \bar J}}_{\mu_2>\nu_K}
\mu_1^{n-|J|-|K|-1}\mu_2^{|K|-1}(\mu_2(d-\nu_J)-\nu_K(\mu_1+\mu_2)).
$$

Type 4. $I=\{2,3\}$, $i_k=3$, $\bar I_1=I_2=\bar I_2=\varnothing$, $I_1=\{2\}$.

By Theorem 2.1, the contribution of such an intersection equals
$$
\multline
\binom{n+1}{|J|+1}(\mu_2+\mu_3-\nu_J)
P_{\nu_{\bar J}-\delta,\delta;\nu(\bar J)}
(\mu_1,\mu_2+\mu_3-\nu_J,\nu(\bar J))\times\\
P_{\mu_{2},\nu_J- \mu_{2}+\delta;\nu(J),\delta}
(\mu_{2},\mu_3,\nu(J),\mu_2+\mu_3-\nu_J).
\endmultline
$$
The first polynomial in the above expression we have already
encountered in the case of two-part double Hurwitz numbers, and its
value is equal to  $(n-|J|)!\mu_1^{n-|J|-1}$. The second polynomial
corresponds to arbitrary two-part
double Hurwitz numbers of total degree $\mu_2+\mu_3$ with $|J|+1$ zeros. By
Corollary~2.2, such numbers are given by
$$
\multline
P_{\mu_{2},\nu_J-\mu_{2}+\delta;\nu(J),\delta}
(\mu_{2},\mu_3,\nu(J),\mu_2+\mu_3-\nu_J)=\\
(|J|+1)!\mathop{\exsum_{K\subseteq J}}_{\nu_J-\mu_2>\nu_K}(\mu_3-\nu_K)
\mu_2^{|J|-|K|}\mu_3^{|K|-1}+\\
(|J|+1)!\mathop{\exsum_{K\subseteq J}}_{\nu_J-\mu_2>\nu_K}(\mu_3-\nu_K)
(\mu_3-\nu_K-(\mu_2+\mu_3-\nu_J))\mu_2^{|J|-|K|-1}\mu_3^{|K|-1};
\endmultline
$$
the first sum in the right hand side
corresponds to resonances not involving the last $y$-coordinate, and
the second one to those involving this coordinate. Therefore, the
total contribution of all intersections of type 4 equals
$$
\multline
(n+1)!\sum_{\mu_2<\nu_J<\mu_2+\mu_3}(\mu_2+\mu_3-\nu_J)\mu_1^{n-|J|-1}
\times\\
\mathop{\exsum_{K\subseteq J}}_{\nu_J-\mu_2>\nu_K}(\mu_3-\nu_K)
\mu_2^{|J|-|K|-1}\mu_3^{|K|-1}(\nu_J\mu_3-\nu_K(\mu_2+\mu_3)).
\endmultline
$$

Collecting all the summands and taking into account the contribution
of the totally negative chamber, we get the following result.

\proclaim{Corollary 2.3} The three-part double Hurwitz numbers are given by
$$
\multline
\frac{h^0_{\mu_1,\mu_2,\mu_3;\nu}}{(n+1)!}=
\exsum_{\mu_2>\nu_J}(\mu_2-\nu_J)\mu_1^{n-|J|-1}\mu_2^{|J|-1}A_J+\\
\sum_{\mu_3>\nu_J}(\mu_3-\nu_J)\mu_3^{|J|-1}
\mathop{\exsum_{K\subseteq \bar J}}_{\mu_2>\nu_K}
\mu_1^{n-|J|-|K|-1}\mu_2^{|K|-1}B_{JK}+\\
\sum_{\mu_2<\nu_J<\mu_2+\mu_3}(\mu_2+\mu_3-\nu_J)\mu_1^{n-|J|-1}
\mathop{\exsum_{K\subseteq J}}_{\nu_J-\mu_2>\nu_K}(\mu_3-\nu_K)
\mu_2^{|J|-|K|-1}\mu_3^{|K|-1}C_{JK},
\endmultline
$$
where $A_J=d\mu_2-\nu_J(\mu_1+\mu_2)$, 
$B_{JK}=(d-\nu_J)\mu_2-\nu_K(\mu_1+\mu_2)$, and 
$C_{JK}=\nu_J\mu_3-\nu_K(\mu_2+\mu_3)$.
\endproclaim

\heading \S 3. Proofs \endheading
\subheading{3.1. Integral representation for double Hurwitz numbers} 
Define the {\it Hurwitz space\/} $\Hu^0_{\mu;\nu}$ as the space of 
degree $d$ meromorphic functions on genus $0$ curves having $m+n-2$ simple
critical values $z_1,\dots,z_{m+n-2}$ and monodromies given by $\mu$
and $\nu$ over two other points $x$ and $y$; the functions are
considered modulo $SL(2,\C)$-action in the image. 
We assume that all
preimages of the points $x, y, z_1,\dots,z_{m+n-2}$ are labelled. The
Lyashko--Looijenga map $\ll$ that takes a function $f\in \Hu^0_{\mu;\nu}$
to the points $x, y, z_1,\dots,z_{m+n-2}$ can be viewed as an
unramified covering of degree $h^0_{\mu;\nu}(d-2)!^{m+n-2}$ over the
moduli space $\Cal M_{0,m+n}$. It is well known that $\ll$ extends
continuously to the mapping $\bll:\bHu^0_{\mu;\nu}\to \MO_{0,m+n}$,
where $\bHu^0_{\mu;\nu}$ is the compactification of $\Hu^0_{\mu;\nu}$
by admissible covers. On the other hand, let $\st:\Hu^0_{\mu;\nu}\to
\Cal M_{0,m+n+(m+n-2)(d-1)}$ be the mapping that takes a function $f$ to
the set of all preimages of the points $x, y, z_1,\dots,z_{m+n-2}$
and let  $\pi:\Cal M_{0,m+n+(m+n-2)(d-1)}\to \Cal M_{0,m+n}$ be the
projection that forgets the preimages of the points
$z_1,\dots,z_{m+n-2}$. Both these mappings extend continuously to the mappings 
between the compactified spaces $\st:\bHu^0_{\mu;\nu}\to
\MO_{0,m+n+(m+n-2)(d-1)}$ and $\pi:\MO_{0,m+n+(m+n-2)(d-1)}\to \MO_{0,m+n}$.
Denote by $x_1,\dots,x_m$ the preimages of $x$
having multiplicities  $\mu_1,\dots,\mu_m$,
and by  $y_1,\dots,y_n$ the preimages of $y$ having multiplicities
$\nu_1,\dots,\nu_n$. Finally, denote by $D$ 
the divisor on $\MO_{0,m+n+(m+n-2)(d-1)}$ whose generic point is a
two-component curve such that $x_1$ lies on one component and
$x_2,\dots,x_m,y_1,\dots,y_n$ lie on the other component. 

\proclaim{Lemma 3.1} One has
$$
\multline
h^0_{\mu;\nu}=
(m+n-2)!\mu_1^{m+n-3}+\\
(d-2)!^{2-m-n}\sum_{u+v=m+n-4}\mu_1\int_{\st(\bHu^0_{\mu;\nu})}D
\left(\mu_1\pi^*\psi(x_1)\right)^{u}
(\st_*\bll^*\psi(x))^{v}.
\endmultline \tag 3.1
$$
\endproclaim

\demo{Proof} Indeed, the degree of $\bll$ equals 
$\int_{\bHu^0_{\mu;\nu}}\bll^*\omega/\int_{\MO_{0,m+n}}\omega$ for any
volume form $\omega$ on $\MO_{0,m+n}$.
It is well known that $\psi(x)^{m+n-3}$ is a volume form 
on $\MO_{0,m+n}$, and that
$\int_{\MO_{0,m+n}}\psi(x)^{m+n-3}=1$. Therefore, 
$$
h^0_{\mu;\nu}(d-2)!^{m+n-2}=
\int_{\bHu^0_{\mu;\nu}}(\bll^*\psi(x))^{m+n-3}.
$$
Besides, by the Ionel lemma
\cite{Io}, $\bll^*\psi(x)=\mu_1\st^*\tilde\psi(x_1)$, where tilde
stands for the classes on $\MO_{0,m+n+(m+n-2)(d-1)}$, and by the
standard pullback formula,
$\tilde\psi(x_1)=\pi^*\psi(x_1)+D$. Therefore, on $\st(\bHu^0_{\mu;\nu})$
one can write
$$
\multline
(\st_*\bll^*\psi(x))^{m+n-3}=\left(\mu_1\pi^*\psi(x_1)\right)^{m+n-3}+\\
\mu_1 D\sum_{u+v=m+n-4}\left(\mu_1\pi^*\psi(x_1)\right)^{u}
(\st_*\bll^*\psi(x))^{v}.
\endmultline 
$$
Finally, $\pi$ is a covering of degree $(m+n-2)!(d-2)!^{m+n-2}$, hence
$$
\multline
\int_{\st(\bHu^0_{\mu;\nu})}\left(\mu_1\pi^*\psi(x_1)\right)^{m+n-3}=\\
(m+n-2)!(d-2)!^{m+n-2}\mu_1^{m+n-3}\int_{\MO_{0,m+n}}\psi(x_1)^{m+n-3}=\\
(m+n-2)!(d-2)!^{m+n-2}\mu_1^{m+n-3},
\endmultline 
$$
and the result follows.
\qed  
\enddemo

\subheading{3.2. Encoding irreducible components of
$\st(\bHu^0_{\mu;\nu})\cap D$} 
In view of Lemma~3.1, we will be interested in the description
of the irreducible components of $\st(\bHu^0_{\mu;\nu})\cap D$.  
Let $f$ be a function whose image belongs to this intersection.
Points $x$ and $y$ on the target curve of $f$
belong to different components. Moreover, the number of components is
exactly two, for the dimensionality reasons. Therefore, the source
curve of  $f$ has a number of double points, which are all mapped to
the unique double point on the target curve. The components of the
source curve are of two types: those covering the component of the
target curve containing $x$, and those covering the component of the
target curve containing $y$. Each component of the first type contains
one or more preimages of $x$, and each component of the second type
contains one or more preimages of  $y$. Finally, the component of the source
curve containing $x_1$ does not contain any other preimages of $x$.

The irreducible components of  $\st(\bHu^0_{\mu;\nu})\cap D$
can be encoded by geometric trees in the
following way. Consider two arbitrary partitions
$[1,m]=\sqcup_{i=1}^{k}I_i$, $[1,n]=\sqcup_{j=1}^{l}J_j$, such that
all parts $I_i$ and $J_j$ are nonempty and $I_1=\{1\}$.
Let $T$ be a tree viewed as a bipartite graph with
the vertices $I_1,\dots,I_k$ in one
part and $J_1,\dots,J_l$ in the other part and let  $\gamma_e$ be a
weight assigned to an edge $e$ of $T$  
in such a way that the sum of $\gamma_e$ over all
edges incident to an arbitrary vertex $I_i$ equals $\mu_{I_i}$, and 
the sum of $\gamma_e$ over all
edges incident to an arbitrary vertex $J_j$ equals
$\nu_{J_j}$. Evidently, $\gamma_e$ are defined by the above condition
in a unique way, and each $\gamma_e$ is an integer. We say that  $T$
is a {\it geometric tree\/} if the following three conditions are satisfied:

(i) at least one among $I_1,\dots,I_k$ is not a singleton;

(ii) all $\gamma_e$ are positive;

(iii) $I_1$ is a leaf of $T$.

The first condition follows from the fact that the target curve of $f$
is stable. A similar condition for $J_1,\dots,J_l$ follows from (iii)
and hence is omitted. The second condition follows from the fact that 
edges of $T$ correspond to double points of the
source curve and weights are the multiplicities at these double points.
The third condition is the definition of $D$. 
The set of all geometric trees is denoted by $\T_{\mu;\nu}$, the
irreducible component of $\st(\bHu^0_{\mu;\nu})\cap D$ encoded by a
geometric tree $T\in \T_{\mu;\nu}$ is denoted by $D_T$.

\remark{Example} Let $\mu=(4,2,1)$, $\nu=(5,2)$. This data defines three
geometric trees, presented on Fig.~1a. Three examples of non-geometric
trees defined by the same data are given on  Fig.~1b. The vertices are
labelled by the corresponding subsets of  $\mu$ and  $\nu$, and the
edges are labelled by the weights  $\gamma_e$.
\endremark

\vskip 15pt
\centerline{\hbox{\epsfysize=2.3cm\epsfbox{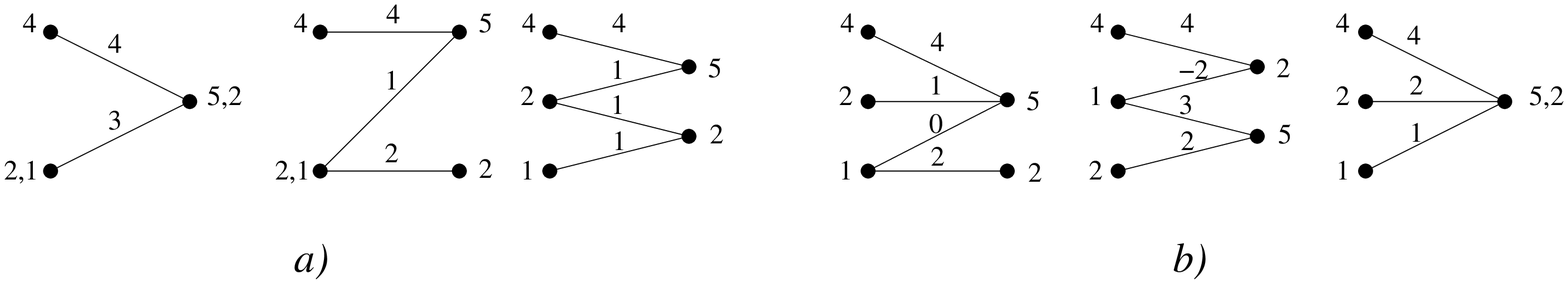}}}
\midspace{0.1mm} \caption{Fig.~1. a) Geometric trees;
b) Several non-geometric trees}

\proclaim{Lemma 3.2} The set of geometric trees is an invariant of a chamber.
\endproclaim

\demo{Proof} It is enough to check that as $(\mu,\nu)$ varies within a
chamber, the weights of the edges of a given graph remain
positive. Any edge $e$ of the initial tree $T$ defines two subtrees of
$T$: the connected components of $T\setminus e$. Denote by $I$, 
$\bar I$ the unions of their $\mu$-vertices, and by $J$, $\bar J$ the
unions of their $\nu$-vertices. It
follows immediately from the definition of weights that the weight of
$e$ can be written as $\mu_I-\nu_J$. Clearly, this quantity does
not change its sign inside a chamber.
\qed
\enddemo

The integral
featuring in~(3.1) can be rewritten as
$$
\multline
\int_{\st(\bHu^0_{\mu;\nu})}D\left(\mu_1\pi^*\psi(x_1)\right)^{u}
(\st_*\bll^*\psi(x))^{v}=\\
\sum_{T\in T_{\mu;\nu}}\delta_T
\int_{D_T}\left(\mu_1\pi^*\psi(x_1)\right)^{u}
(\st_*\bll^*\psi(x))^{v},
\endmultline \tag 3.2
$$
where $\delta_T$ is the multiplicity arising from the non-transversal
intersection of $\st(\bHu^0_{\mu;\nu})$ and $D$. This multiplicity can
be calculated as follows.

 \proclaim{Lemma 3.3} 
The multiplicity $\delta_T$ is given by
$$
\delta_T=\frac1{\mu_1}\prod_{e\in T}\gamma_e.
$$
\endproclaim

\demo{Proof} Let $f\in\bHu^0_{\mu;\nu}$ be a function such that 
$\st(f)\in D_T$. To find $\delta_T$, we have to parametrize the local
universal deformation space in a small neighborhood of
$C=\st(f)$. Recall that $C_0=\ll(f)$ is a two-component curve. Following
\cite{HM, pp.~61,62}, we choose coordinates $u, v$ at the double point
of $C_0$ so that $C_0$ is given locally by $uv=0$. Then, for $t$
suitably chosen, the universal deformation of $C_0$ is given
locally near the double point by the equation $uv=t$. Similarly, at each double
point $p_e$, $e\in T$, of $C$ the universal deformation of 
$\st(\bHu^0_{\mu;\nu})$ is defined by $x_ey_e=t_e$ with
$u=x_e^{\gamma_e}$, $v=y_e^{\gamma_e}$. We get therefore
$t_e^{\gamma_e}=t$, $e\in T$. Recall that $D_T$ locally is given by
the equation $t_{e_1}=0$, where $e_1$ is the edge incident to the leaf
$I_1$. Therefore, $\delta_T$ is equal to the multiplicity of the point
$\{t_e=0, e\in T, e\ne e_1\}$ in the curve $\{t_e^{\gamma_e}=t, e\in
T, e\ne e_1\}$. The latter is equal to $\frac1{\mu_1}\prod_{e\in
T}\gamma_e$, since $\gamma_{e_1}=\mu_1$.   
\qed
\enddemo

\subheading{3.3. Essential geometric trees} We say that a geometric
tree $T\in \T_{\mu;\nu}$ is {\it essential\/} if all vertices $J_j$
except for the one connected to $I_1$ are singletons. The set of all
essential geometric trees for a given pair $(\mu,\nu)$ is denoted
$\E_{\mu;\nu}$. As follows immediately from Lemma~3.2, $\E_{\mu;\nu}$
is an invariant of a chamber.
The importance of
this notion is revealed in the following statement. Denote by
$\I(T,u,v)$ the integral in the right hand side of~(3.2).

 \proclaim{Lemma 3.4} For any inessential tree  $T\in
 \T_{\mu;\nu}\setminus\E_{\mu;\nu}$ 
 one has $\I(T,u,v)=0$.
\endproclaim

\demo{Proof} For any vertex $I_i$ of $T$ denote by  $\gamma(I_i)$ the 
partition of $\mu_{I_i}$ formed by the weights of the edges incident
to $I_i$, and by  $\deg I_i$ the number of such edges; $\gamma(J_j)$
and $\deg J_j$ will have a similar meaning for the vertex $J_j$.
Let $\bHu^0_{I_1,\dots,I_k}$ denote the compactified Hurwitz
space of degree $d$ meromorphic functions on disconnected genus~$0$
curves with $k$ connected components such that on the $i$th component
the function has $|I_i|+\deg I_i-2$ simple critical values and
monodromies given by $\mu(I_i)$ and $\gamma(I_i)$ over the points $x$
and $y$. Note that this space does not coincide with the direct
product of the spaces $\bHu^0_{\mu(I_i);\gamma(I_i)}$; they both are
obtained from the same space by taking quotient by two different
actions: the first one by the action of $\C^*$, and the second one by
the action of the direct product of several copies of $\C^*$. 
Consider a natural mapping 
$$
\rho_T:D_T\to \prod_{j=1}^l
\bHu^0_{\nu(J_j);\gamma(J_j)}\times\bHu^0_{I_1,\dots,I_k}. \tag 3.3
$$
Denote by $\bHu$ the space in the right hand side of the above formula.
It is easy to see that the class to be integrated over $D_T$ is a
pullback under $\rho_T$ of a certain class on $\bHu$. 
Therefore, $\I(T,u,v)$ will vanish whenever the dimenson of $\bHu$
differs from that of  $D_T$.

Let us calculate the dimension of $\bHu$. Clearly,
$\dim\bHu^0_{\nu(J_j);\gamma(J_j)}=\max\{|J_j|+\deg
 J_j-3,0\}$. Since $|J_j|+\deg J_j-3$ equals $-1$ if $J_j$ is a singleton
and is nonnegative otherwise, we see that
the total dimension of the direct product of
all $\bHu^0_{\nu(J_j);\gamma(J_j)}$ equals 
$$
n+(k+l-1)-3l+l^*=n+k-2l+l^*-1,
$$
where $l^*$ is the number of singletons among $J_1,\dots,J_l$. 
Further, 
$$
\dim\bHu^0_{I_1,\dots,I_k}=\sum_{i=1}^k(|I_i|+\deg I_i-2)-1=
m+(k+l-1)-2k-1=m+l-k-2.
$$
Finally, $\dim D_T=m+n-4$. Equating the two dimensions we get
$l-l^*=1$, which means that $T$ is essential.
\qed
\enddemo

The value of the integral  $\I(T,u,v)=0$ for essential trees is
calculated in the following proposition.

 \proclaim{Lemma 3.5} For any essential tree  $T\in \E_{\mu;\nu}$
 one has
$$
\I(T,u,v)=\cases (d-2)!^{r}
r!\mu_1^u{\dsize\prod_{i=1}^k\frac
{h^0_{\mu(I_i);\gamma(I_i)}}{r_i!}\prod_{j\in \bar J_1}\nu_{j}^{-1}}\qquad
&\text{if $u=|J_1|+\deg J_1-3$}\\
0\qquad &\text{otherwise},
\endcases
$$
where $r=m+n-2$ and $r_i=|I_i|+\deg I_i-2$.
\endproclaim

\demo{Proof} Define $\rho_T$ as in~(3.3). Since $T$ is essential,
$\rho_T$ is a finite map, and its degree is equal to
$\binom{m+n-2}{u+1}[(d-2)!/(\nu_{J_1}-2)!]^{u+1}$. 
Using a Fubini-type argument we get
$$
\multline
\I(T,u,v)=\binom{m+n-2}{u+1}\left(\frac{(d-2)!}
{(\nu_{J_1}-2)!}\right)^{u+1}\times\\
\int_{\bHu^0_{\nu(J_1);\gamma(J_1)}}\left(\mu_1\pi^*\psi(x_1)\right)^{u}
\prod_{j=2}^l\int_{\bHu^0_{\nu(J_j);\gamma(J_j)}}1
\int_{\bHu^0_{I_1,\dots,I_k}}(\tilde\ll^*\psi(x))^{v},
\endmultline \tag 3.4
$$
where $\tilde\ll$ is the Lyashko--Looijenga map for the space
$\bHu^0_{I_1,\dots,I_k}$. The first integral in~(3.4)
vanishes whenever $u$ is distinct from the dimension of
$\bHu^0_{\nu(J_1);\gamma(J_1)}$, which is equal to $|J_1|+\deg
J_1-3$. 

From now on we assume that $u=|J_1|+\deg J_1-3$. Than the first
integral in~(3.4) equals 
$\mu_1^u\int_{\MO_{0,|J_1|+\deg J_1}}\psi(x_1)^{u}$ times the
transversal multiplicity of the projection $\pi$ at a point of $D_T$.
The latter is equal to $(u+1)!(\nu_{J_1}-2)!^{u+1}$, and the integral over
the moduli space equals~$1$, so the first integral in~(3.4) equals 
$\mu_1^u(u+1)!(\nu_{J_1}-2)!^{u+1}$.  

The dimension of
the space $\bHu^0_{I_1,\dots,I_k}$ equals $v$, and hence the
last integral in~(3.4) gives the 
degree of $\tilde\ll$, similarly to the argument in Section~3.1.  
The degree of $\tilde\ll$ equals $(d-2)!^{v+1}\binom{v+1}{r_1
\cdots r_k}$ times the corresponding
double Hurwitz number of disconnected coverings. The latter is equal to the
product of double Hurwitz numbers over the components of the covering.
Finally, since $T$ is essential, all $J_j$ for $j=2,\dots,l$ are
singletons, hence the $j$th among the $l-1$ intermediate
integrals in~(3.4) equals $1/\nu_{J_j}$,  and their product is 
$\prod_{j\in \bar J_1}\nu_{j}^{-1}$. Multiplying all the above factors and
taking into account that $u+v=m+n-4$, we 
arrive to the desired result.
\qed
\enddemo

Let us take a more precise look at essential geometric trees. The
structure of such trees is very simple, as shown on Fig.~2. We denote
$\mu$-vertices by circles and $\nu$-vertices by squares. Singletons
are white and non-singletons are black. An essential geometric tree
has a unique black square vertex denoted $J_1$. 
All the other black vertices are circles; we denote them
$\widehat{I}_i$, $i=1,\dots,\widehat{k}$. 
Note that the collection of all $\widehat{I}_i$'s forms
a proper subset of the initial collection of all ${I}_i$'s. For any
black circle vertex $\widehat{I}_i$ we denote by $\widehat{J}_i$ the
(possibly empty) union of all white squares incident to it. Note that
in general $\widehat{J}_i$ does not coincide with any of the initial ${J}_j$.

\vskip 15pt
\centerline{\hbox{\epsfysize=3cm\epsfbox{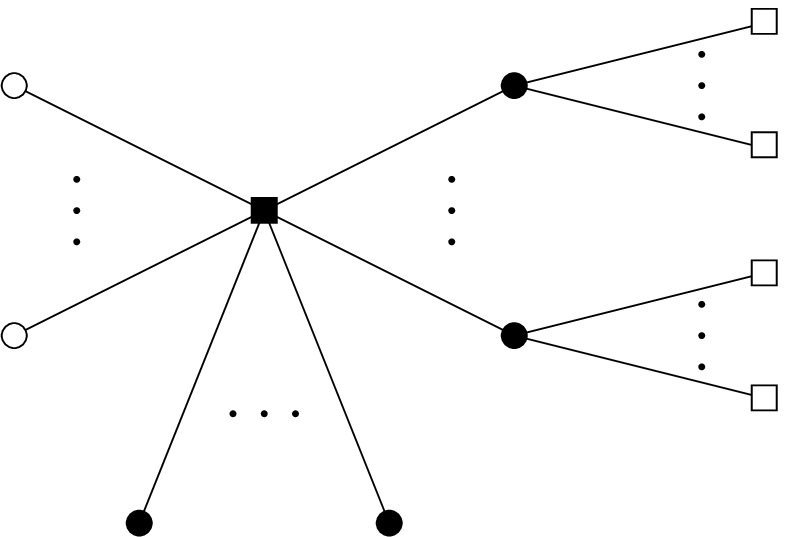}}}
\midspace{0.1mm} \caption{Fig.~2. Essential geometric trees}

We can now
represent double Huritz numbers as a sum over the set of essential
geometric trees of products of double Hurwitz numbers of a smaller
size. For unification purposes, we extend the set of essential
geometric trees by adding the tree with $m$ white circle vertices
corresponding to $\mu_i$, and no black circle or white square
vertices. The extended set is denoted $\E^*_{\mu;\nu}$.

\proclaim{Theorem 3.6} Double Hurwitz numbers are given by
$$
h^0_{\mu;\nu}=(m+n-2)!   
\sum_{T\in \E^*_{\mu;\nu}}
\mu_1^{u_T}\prod_{i=1}^{\widehat{k}}(\mu_{\widehat{I}_i}-\nu_{\widehat{J}_i})
\frac{h^0_{\mu(\widehat{I}_i);\nu^*(\widehat{J}_i)}}
{(|\widehat{I}_i|+|\widehat{J}_i|-1)!},
\tag 3.5
$$
where $u_T=|J_1|+\deg J_1-3$ and $\nu^*(\widehat{J}_i)$ is obtained from 
$\nu(\widehat{J}_i)$ by insertion of
$\mu_{\widehat{I}_i}-\nu_{\widehat{J}_i}$ at the proper place.
\endproclaim
 
\demo{Proof} Follows from Lemmas 3.1, 3.3--3.5 once we note that edges
connecting $J_1$ with white circle vertices contribute $\mu_i$ to
$\delta_T$ and $\mu_i^{-1}$ to the product of double Hurwitz numbers in the
formula for $\I(T,u,v)$ in Lemma~3.5, so that their contributions are
cancelled. Further, edges
connecting $J_1$ with black circle vertices contribute
$\mu_{\widehat{I}_i}-\nu_{\widehat{J}_i}$ to
$\delta_T$. Finally, edges connecting black circle vertices to white
square vertices contribute $\nu_j$ to $\delta_T$, and so their
contribution is cancelled with the last factor in the
formula for $\I(T,u,v)$ in Lemma~3.5. The added tree accounts for the
first term in the right hand side of~(3.1), since in this case
$|J_1|=n$ and $\deg J_1=m$.
\qed
\enddemo

\subheading{3.4. Proofs of Theorems~1.1--1.3} Let us start from the
following observation.

\proclaim{Proposition 3.7} Resonances for the double Hurwitz numbers 
$h^0_{\mu(\widehat{I}_i);\nu^*(\widehat{J}_i)}$ correspond bijectively to
resonances $x_I=y_J$ for  $h^0_{\mu;\nu}$ with $I\subset
\widehat{I}_i$, $J\subset \widehat{J}_i$. 
\endproclaim

The proof of Theorem~1.1 follows immediately from Lemma~3.2,
Theorem~3.6 and Proposition~3.7 by induction over $m+n$. The base of
induction is formed by the cases when either $m$ or $n$ equals~$1$,
and there are no resonances.

To prove Theorem~1.2, observe that by Proposition~3.7, (3.5) can be
rewritten as
$$
P_{\mu;\nu}=(m+n-2)!  
\sum_{T\in \E^*_{\mu;\nu}}
x_1^{u_T}\prod_{i=1}^{\widehat{k}}(x_{\widehat{I}_i}-y_{\widehat{J}_i})
\frac{P_{\mu(\widehat{I}_i);\nu^*(\widehat{J}_i)}}
{(|\widehat{I}_i|+|\widehat{J}_i|-1)!}.
\tag 3.6
$$
We denote
by $Q^T_{\mu;\nu}$ the contribution of a tree $T$ to the right hand
side of~(3.6).

Assume that $(\mu,\nu)\in C$ and $(\mu',\nu')\in C'$. Clearly,
$\E^*_{\mu';\nu'}\subset \E^*_{\mu;\nu}$.
Therefore, the expression for
$P_{\mu;\nu}$ may change for two reasons: 

1) birth of new essential geometric trees, and 

2) changes in the expressions for $Q^T_{\mu;\nu}$.

Let $\E_{res}$ be the set of essential geometric trees
$T\in\E_{\mu';\nu'}$ such that $Q^T_{\mu;\nu}\ne Q^T_{\mu';\nu'}$, and
let $\E_{new}=\E_{\mu;\nu}\setminus\E_{\mu';\nu'}$. 

\proclaim{Lemma 3.8} There exists a bijection between
$\E_{res}\cup\E_{new}$ and $\E^*_{\mu(\bar I),\mu_I-\nu_J;\nu(\bar
J)}$.
\endproclaim

\demo{Proof} The bijection $\Phi$ is presented on Figure~3 below. The upper
part of the figure corresponds to the trees in $\E_{new}$, the lower part
corresponds to the trees in $\E_{res}$. The asterisque stands for the
additional index corresponding to $\mu_I-\nu_J$ in  $\E^*_{\mu(\bar
I),\mu_I-\nu_J;\nu(\bar J)}$. 

\vskip 15pt
\centerline{\hbox{\epsfxsize=12cm\epsfbox{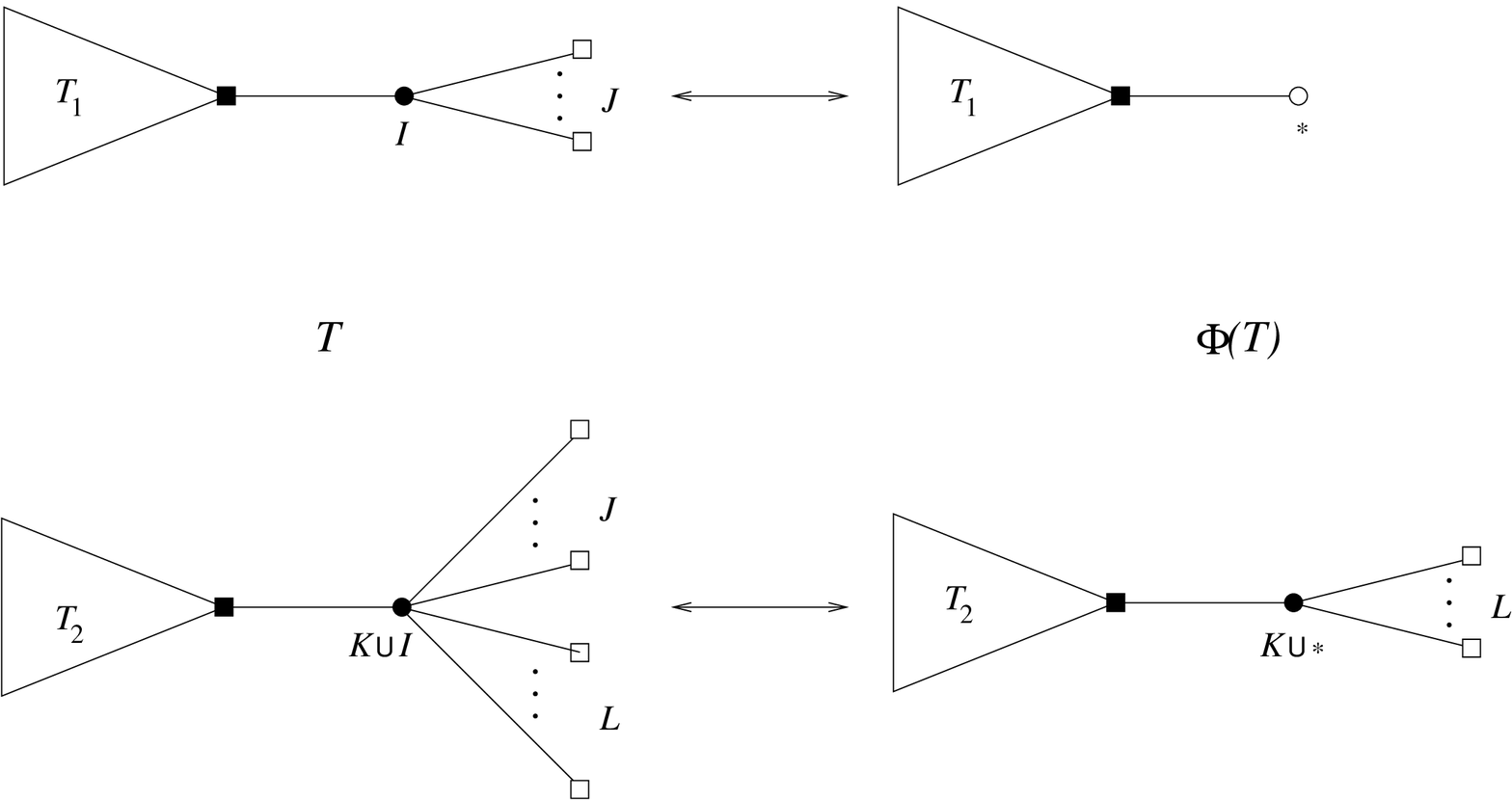}}}
\midspace{0.1mm} \caption{Fig.~3. The definition of the bijection $\Phi$}

\qed
\enddemo
 
For a tree $T\in\E_{new}$ we have
$$
Q^T_{\mu;\nu}=
(m+n-2)!
x_1
\frac{Q^{T_1}_{\mu(\bar I);\nu(\bar J)}}{(m-|I|+n-|J|-2)!}
(x_I-y_J)\frac{P_{\mu(I);\nu(J),\mu_I-\nu_J}}{(|I|+|J|-1)!}.
$$
Besides, it follows from the definition of $\Phi$ that
$$
\frac{x_1Q^{T_1}_{\mu(\bar I);\nu(\bar J)}}{(m-|I|+n-|J|-2)!}=
\frac{Q^{\Phi(T)}_{\mu(\bar I),\mu_I-\nu_J;\nu(\bar
J)}}{(m-|I|+n-|J|-1)!},
$$
and hence
$$
Q^T_{\mu;\nu}=\binom{m+n-2}{|I|+|J|-1}Q^{\Phi(T)}_{\mu(\bar
I),\mu_I-\nu_J;\nu(\bar J)}(x_I-y_J)P_{\mu(I);\nu(J),\mu_I-\nu_J}.
$$

For a tree $T\in\E_{res}$ we have
$$
\multline
Q^T_{\mu;\nu}-Q^T_{\mu';\nu'}=(m+n-2)!x_1
\frac{Q^{T_2}_{\mu(\overline{I\cup K});\nu(\overline{J\cup L})}}
{(m-|I|-|K|+n-|J|-|L|-2)!}\times\\
(x_I+x_K-y_J-y_L)(P_{\mu(I\cup K);\nu^*(J\cup L)}- P_{\mu'(I\cup
K);{\nu'}^*(J\cup L)}).
\endmultline
$$ 
By induction over $m+n$, the latter difference equals
$$
\binom{|I|+|K|+|J|+|L|-1}{|I|+|J|-1}(x_I-y_J)P_{\mu(K),\mu_I-\nu_J;
\nu^*(L)}P_{\mu(I);\nu(J),\mu_I-\nu_J}.
$$
Besides, it follows from the definition of $\Phi$ that
$$
\multline
\frac{x_1Q^{T_2}_{\mu(\overline{I\cup K});\nu(\overline{J\cup L})}}
{(m-|I|-|K|+n-|J|-|L|-2)!}(x_I+x_K-y_J-y_L)
\frac{P_{\mu(K),\mu_I-\nu_J;\nu^*(L)}}{(|K|+|L|)!}=\\
\frac{Q^{\Phi(T)_{\mu(\bar
I),\mu_I-\nu_J;\nu(\bar J)}}}{(m-|I|+n-|J|-1)!},
\endmultline
$$
and hence
$$
Q^T_{\mu;\nu}-Q^T_{\mu';\nu'}=\binom{m+n-2}{|I|+|J|-1}Q^{\Phi(T)}_{\mu(\bar
I),\mu_I-\nu_J;\nu(\bar J)}(x_I-y_J)P_{\mu(I);\nu(J),\mu_I-\nu_J}.
$$
Summing up the contributions of all trees in $\E_{new}\cup\E_{res}$
and using Lemma~3.8 we get
$$
\multline
P_{C}-P_{C'}=\\
\binom{m+n-2}{|I|+|J|-1}(x_I-y_J)
P_{\mu(I);\nu(J),\mu_I-\nu_J}
\sum_{T\in\E^*_{\mu(\bar I),\mu_I-\nu_J,\nu(\bar J)}}Q^{\Phi(T)}_{\mu(\bar
I),\mu_I-\nu_J;\nu(\bar J)}=\\
\binom{m+n-2}{|I|+|J|-1}(x_I-y_J)
P_{\mu(I);\nu(J),\mu_I-\nu_J}P_{\mu(\bar I),\mu_I-\nu_J;\nu(\bar J)},
\endmultline
$$
as required. 

It remains to prove Theorem~1.3. We start from expression~(3.5). 
Note that essential geometric trees
for the totally negative chamber have a very simple structure: they do
not have white square vertices. Therefore, $\mu_{\widehat{I}_i}-
\nu_{\widehat{J}_i}= \mu_{\widehat{I}_i}$ and
$h^0_{\mu(\widehat{I}_i); \nu^*(\widehat{J}_i)}=(|\widehat{I}_i|-1)!
\mu_{\widehat{I}_i}^{|\widehat{I}_i|-2}$. Finally, $|J_1|=n$, and
hence the double Hurwitz numbers for the totally negative chamber are
given by
$$
(m+n-2)!\mu_1^{n-1}\sum_{k=1}^{m-1}\mu_1^{k-1}
\sum_{I_1,\dots,I_k}\prod_{i=1}^k\mu_{I_i}^{|I_i|-1},
$$
where the inner sum is taken over all unordered partitions of $[2,m]$ into 
$k$ nonempty parts $I_1, \dots, I_k$. It is easy to see that
$$
\sum_{k=1}^{m-1}x_1^k\sum_{I_1,\dots,I_k}\prod_{i=1}^k x_{I_i}^{|I_i|-1}=
x_1(x_1+\cdots+x_m)^{m-2}, \tag 3.7
$$
since both parts of the above formula enumerate trees on
$[1,n]$ rooted at~$1$ classified according to the degrees 
of the vertices (see e.g. \cite{Pi, Cayley's tree volume formula}. 
Therefore, the double Hurwitz numbers in
question equal $(m+n-2)!\mu_1^{n-1}(\mu_1+\cdots+\mu_m)^{m-2}$, 
as required.

\remark{Remark} Identity~(3.7) is a Hurwitz type multinomial identity, see
\cite{Pi} and references therein. 
Identities of this kind were discovered by Hurwitz in \cite{Hu} and, 
apparently, used by him in his studies of Hurwitz numbers (see \cite{St} 
for a conjectural reconstruction of the original Hurwitz derivation for
$h^0_\mu$). It is interesting to note that such identities arose again 
recently in connection with Gromov--Witten invariants in \cite{GOP}.
\endremark


\Refs
\widestnumber \key{ELSV2}

\ref\key{ELSV1}
\by T.~Ekedahl, S.~Lando, M.~Shapiro, and A.~Vainshtein
\paper On Hurwitz numbers and Hodge integrals
\jour C. R. Acad. Sci. Paris S\'er. I Math.
\vol 328
\yr 1999
\pages   1175-1180
\endref

\ref\key{ELSV2}
\by T.~Ekedahl, S.~Lando, M.~Shapiro, and  A.~Vainshtein
\paper Hurwitz numbers and intersections on moduli spaces of curves
\jour Invent. Math. 
\vol 146
\yr 2001
\pages  297--327
\endref

\ref\key{GOP}
\by E.~Getzler, A.~Okounkov, and R.~Pandharipande
\paper Multipoint series of Gromov--Witten invariants of $\CP^1$
\jour Lett. Math. Phys. 
\vol 62
\yr 2002
\pages  159--170
\endref

\ref\key{GJ}
\by I.~Goulden and D.~Jackson
\paper The combinatorial relationship between trees, cacti and
certain connection coefficients for the symmetric group
\jour Europ. J. Comb.
\yr 1992
 \vol 13 
\pages 357--365
\endref


\ref\key{GJVa}
\by I.~Goulden, D.~Jackson, and R.~Vakil
\paper Towards the geometry of double Hurwitz numbers
\jour Adv. Math.
\yr 2005
\vol 198
\pages 43--92 
\endref

\ref\key{GrVa}
\by T.~Graber and R.~Vakil
\paper Hodge integrals and Hurwitz numbers via virtual localization 
\jour Compositio Math. 
\vol 135
\yr 2003
\pages  25--36
\endref

\ref\key{HM}
\by J.~Harris and D.~Mumford
\paper On the Kodaira dimension of the moduli space of curves
\jour Inv. Math. 
\vol 67
\yr 1982
\pages  23--86
\endref


\ref\key{Hu}
\by A.~Hurwitz
\paper \"Uber Abel's Verallgemeinerung der binomischen Formel 
\jour Acta Math
\vol 26
\yr 1902
\pages  199--203
\endref

\ref\key{Io}
\by E.-N.~Ionel
\paper Topological recursive relations in $H^{2g}(\Cal M_{g,n})$
\jour Invent. Math.   
\vol 148
\yr 2002
\pages  627--658
\endref

\ref\key{KL}
\by M.~Kazaryan, S.~Lando
\paper On intersection theory on Hurwitz spaces
\jour Izv. Math.   
\vol 68
\yr 2004
\pages  935--964
\endref

\ref\key{Ku}
\by S.~Kuleshov
\paper Ramified coverings of $S^2$ with two degenerate branching
points
\yr 2003
\finalinfo preprint MPIM2001-62
\endref


\ref\key{Ok}
\by A.~Okounkov
\paper Toda equations for Hurwitz numbers
\jour Math. Res. Lett.
\issue 7
\yr 2000
\pages  447--453
\endref

\ref\key{Pi}
\by J.~Pitman
\paper Forest volume decompositions and Abel--Cayley--Hurwitz multinomial 
expansions
\jour J. Combin. Theory Ser. A
\vol 98
\yr 2002
\pages  175--191
\endref

\ref\key{Sh1}
\by S.~Shadrin
\paper Geometry of meromorphic functions and intersections on moduli 
spaces of curves
\jour Int. Math. Res. Not.
\issue 38
\yr 2003
\pages  2051--2094
\endref

\ref\key{Sh2}
\by S.~Shadrin
\paper Hurwitz numbers of generalized polynomials, and cycles of two-point 
branchings
\jour Uspekhi Mat. Nauk
\vol 58
\yr 2003
\pages  197--198
\endref

\ref\key{Sh3}
\by S.~Shadrin
\paper Polynomial Hurwitz numbers and intersections on $\MO_{0,k}$
\jour Funct. Anal. Appl.
\vol 37
\yr 2003
\pages  78--80
\endref

\ref\key{SSV}
\by B.~Shapiro, M.~Shapiro, and A.~Vainshtein 
\paper Ramified coverings of $S^2$ with one degenerate branching
point and enumeration of edge-oriented graphs
\jour Adv. in Math. Sci.
\vol 34
\yr 1997
\pages 219--228
\endref


\ref\key{St}
\by V.~Strehl 
\paper Minimal transitive products of transpositions: the reconstruction of
a proof of A.~Hurwitz
\jour Sem. Loth. Comb.
\vol 37
\yr 1996
\finalinfo B37c
\endref

\ref\key{Z1}
\by D.~Zvonkine 
\paper Multiplicities of the Lyashko--Looijenga map on its strata 
\jour C. R. Acad. Sci. Paris S\'er. I Math.
\vol 324
\yr 1997
\pages 1349--1353
\endref

\ref\key{Z2}
\by D.~Zvonkine 
\paper Counting ramified coverings and intersection theory on Hurwitz spaces
II (Local structure of Hurwitz spaces and combinatorial results)
\jour Moscow Math. J. 
\yr 2006
\finalinfo to appear; see also math.AG/0304251
\endref



\ref\key{ZL}
\by D.~Zvonkin and S.~Lando
\paper On multiplicities of the Lyashko--Looijenga mapping on strata of the 
discriminant
\jour Funct. Anal. Appl.
\vol 33
\yr 1999
\pages 178--188
\endref

\endRefs
\enddocument